# Darboux Approach to Bertrand Surface Offsets

**Mehmet ÖNDER**
*Celal Bayar University, Faculty of Science and Arts, Department of Mathematics, Muradiye Manisa, Turkey. E-mail:* mehmet.onder@bayar.edu.tr

**Abstract**
In this paper, we study Bertrand surface offsets by considering the dual geodesic trihedron(dual Darboux frame) of the ruled surfaces. We obtain the relationships between the invariants of Bertrand trajectory ruled surfaces. Furthermore, we obtain the conditions for these surface offset to be developable.

## 1. Introduction

In the space, a continuously moving of a straight line generates a special surface called ruled surface. A ruled surface can always be described (at least locally) as the set of points swept by a moving straight line. These surfaces are one of the most fascinating topics of the surface theory and also used in many areas of science such as Computer Aided Geometric Design(CAGD), mathematical physics, moving geometry, kinematics for modeling the problems and model-based manufacturing of mechanical products. Especially, the offsets of the ruled surface have an important role in (CAGD)[6,7]. In [8], Ravani and Ku defined and given the generalization of the theory of Bertrand curve for Bertrand trajectory ruled surfaces on the line geometry. By considering the E. Study mapping, Küçük and Gürsoy have studied the integral invariants of closed Bertrand trajectory ruled surfaces in dual space[5]. They have given some characterizations of Bertrand offsets of trajectory ruled surfaces in terms of integral invariants (such as the angle of pitch and the pitch) of closed ruled surfaces and they have obtained the relationships between the projection areas for the spherical images of Bertrand surface offsets and their integral invariants. They have also found that if the Bertrand offsets of trajectory ruled surfaces are developable then their striction lines are Bertrand partner curves.

In this paper, we examine the Bertrand offsets of the ruled surfaces in view of their dual geodesic trihedrons(dual Darboux frames). Using the dual representation of the ruled surfaces, we give some theorems and the results characterizing Bertrand surface offsets.

## 2. Dual Numbers and Dual Vectors

Let $D = IR \times IR = \{\bar{a} = (a, a^*) : a, a^* \in IR\}$ be the set of the pairs $(a, a^*)$. For $\bar{a} = (a, a^*)$, $\bar{b} = (b, b^*) \in D$ the following operations are defined on $D$:

Equality: $\bar{a} = \bar{b} \Leftrightarrow a = b, \ a^* = b^*$
Addition: $\bar{a} + \bar{b} \Leftrightarrow (a+b, \ a^* + b^*)$
Multiplication: $\bar{a}\bar{b} \Leftrightarrow (ab, \ ab^* + a^*b)$

The element $\varepsilon = (0,1) \in D$ satisfies the relationships

$$\varepsilon \neq 0, \quad \varepsilon^2 = 0, \quad \varepsilon 1 = 1\varepsilon = \varepsilon \tag{1}$$

Let consider the element $\bar{a} \in D$ of the form $\bar{a} = (a, 0)$. Then the mapping $f : D \to IR$, $f(a, 0) = a$ is a isomorphism. So, we can write $a = (a, 0)$. By the multiplication rule we have that

$$\begin{aligned}
\bar{a} &= (a, a^*) \\
&= (a, 0) + (0, a^*) \\
&= (a, 0) + (0, 1)(a^*, 0) \\
&= a + \varepsilon a^*
\end{aligned}$$

Then $\bar{a} = a + \varepsilon a^*$ is called dual number and $\varepsilon$ is called dual unit. Thus the set of dual numbers is given by

$$D = \{\bar{a} = a + \varepsilon a^* : a, a^* \in IR, \varepsilon^2 = 0\} \tag{2}$$

The set $D$ forms a commutative group under addition. The associative laws hold for multiplication. Dual numbers are distributive and form a ring over the real number field[2].

Dual function of dual number presents a mapping of a dual numbers space on itself. Properties of dual functions were thoroughly investigated by Dimentberg[2]. He derived the general expression for dual analytic (differentiable) function as follows

$$f(\bar{x}) = f(x + \varepsilon x^*) = f(x) + \varepsilon x^* f'(x), \tag{3}$$

where $f'(x)$ is derivative of $f(x)$ and $x, x^* \in IR$. This definition allows us to write the dual forms of some well-known functions as follows

$$\begin{cases} \cos(\bar{x}) = \cos(x + \varepsilon x^*) = \cos(x) - \varepsilon x^* \sin(x), \\ \sin(\bar{x}) = \sin(x + \varepsilon x^*) = \sin(x) + \varepsilon x^* \cos(x), \\ \sqrt{\bar{x}} = \sqrt{x + \varepsilon x^*} = \sqrt{x} + \varepsilon \frac{x^*}{2\sqrt{x}}, \quad (x > 0). \end{cases} \tag{4}$$

Let $D^3 = D \times D \times D$ be the set of all triples of dual numbers, i.e.,

$$D^3 = \{\tilde{a} = (\bar{a}_1, \bar{a}_2, \bar{a}_3) : \bar{a}_i \in D, i = 1, 2, 3\}, \tag{5}$$

Then the set $D^3$ is called dual space. The elements of $D^3$ are called dual vectors. Similar to the dual numbers, a dual vector $\tilde{a}$ may be expressed in the form $\tilde{a} = \vec{a} + \varepsilon \vec{a}^* = (\vec{a}, \vec{a}^*)$, where $\vec{a}$ and $\vec{a}^*$ are the vectors of $IR^3$. Then for any vectors $\tilde{a} = \vec{a} + \varepsilon \vec{a}^*$ and $\tilde{b} = \vec{b} + \varepsilon \vec{b}^*$ of $D^3$, the scalar product and the vector product are defined by

$$\langle \tilde{a}, \tilde{b} \rangle = \langle \vec{a}, \vec{b} \rangle + \varepsilon \left( \langle \vec{a}, \vec{b}^* \rangle + \langle \vec{a}^*, \vec{b} \rangle \right), \tag{6}$$

and

$$\tilde{a} \times \tilde{b} = \vec{a} \times \vec{b} + \varepsilon \left( \vec{a} \times \vec{b}^* + \vec{a}^* \times \vec{b} \right), \tag{7}$$

respectively, where $\langle \vec{a}, \vec{b} \rangle$ and $\vec{a} \times \vec{b}$ are the inner product and the vector product of the vectors $\vec{a}$ and $\vec{a}^*$ in $IR^3$, respectively.

The norm of a dual vector $\tilde{a}$ is given by

$$\|\tilde{a}\| = \|\vec{a}\| + \varepsilon \frac{\langle \vec{a}, \vec{a}^* \rangle}{\|\vec{a}\|}, \quad (\vec{a} \neq 0). \tag{8}$$

A dual vector $\tilde{a}$ with norm $1 + \varepsilon 0$ is called dual unit vector. The set of dual unit vectors is given by

$$\tilde{S}^2 = \{\tilde{a} = (a_1, a_2, a_3) \in D^3 : \langle \tilde{a}, \tilde{a} \rangle = 1 + \varepsilon 0\}, \tag{9}$$

and called dual unit sphere[1,3].

## 3. Dual Representation of Ruled Surfaces

In the Euclidean 3-space $IR^3$, an oriented line $L$ is determined by a point $p \in L$ and a unit vector $\vec{a}$. Then, one can define $\vec{a}^* = \vec{p} \times \vec{a}$ which is called moment vector. The value of $\vec{a}^*$ does not depend on the point $p$, because any other point $q$ in $L$ can be given by $\vec{q} = \vec{p} + \lambda \vec{a}$ and then $\vec{a}^* = \vec{p} \times \vec{a} = \vec{q} \times \vec{a}$. Reciprocally, when such a pair $(\vec{a}, \vec{a}^*)$ is given, one recovers the line $L$ as $L = \{(\vec{a} \times \vec{a}^*) + \lambda \vec{a} : \vec{a}, \vec{a}^* \in E^3, \lambda \in IR\}$, written in parametric equations. The vectors $\vec{a}$ and $\vec{a}^*$ are not independent of one another and they satisfy the following relationships

$$\langle \vec{a}, \vec{a} \rangle = 1, \quad \langle \vec{a}, \vec{a}^* \rangle = 0 \tag{10}$$

The components $a_i$, $a_i^*$ ($1 \leq i \leq 3$) of the vectors $\vec{a}$ and $\vec{a}^*$ are called the normalized Plucker coordinates of the line $L$. From (6), (9) and (10) we see that the dual unit vector $\tilde{a} = \vec{a} + \varepsilon \vec{a}^*$ corresponds to the line $L$. This correspondence is known as E. Study Mapping: There exists a one-to-one correspondence between the vectors of dual unit sphere $\tilde{S}^2$ and the directed lines of the space $IR^3$. By the aid of this correspondence, the properties of the spatial motion of a line can be derived. Hence, the geometry of ruled surface is represented by the geometry of dual curves lying on the dual unit sphere $\tilde{S}^2$.

The angle $\bar{\theta} = \theta + \varepsilon \theta^*$ between two dual unit vectors $\tilde{a}, \tilde{b}$ is called *dual angle* and defined by

$$\langle \tilde{a}, \tilde{b} \rangle = \cos \bar{\theta} = \cos \theta - \varepsilon \theta^* \sin \theta . \tag{11}$$

By considering the E. Study Mapping, the geometric interpretation of dual angle is that, $\theta$ is the real angle between the lines $L_1, L_2$ corresponding to the dual unit vectors $\tilde{a}, \tilde{b}$, respectively, and $\theta^*$ is the shortest distance between those lines[3].

Let now $(\tilde{k})$ be a dual curve represented by the dual vector $\tilde{e}(u) = \vec{e}(u) + \varepsilon \vec{e}^*(u)$. The unit vector $\vec{e}$ draws a curve on the real unit sphere $S^2$ and is called the (real) indicatrix of $(\tilde{k})$. We suppose throughout that it is not a single point. We take the parameter $u$ as the arc-length parameter $s$ of the real indicatrix and denote the differentiation with respect to $s$ by primes. Then we have $\langle \vec{e}', \vec{e}' \rangle = 1$. The vector $\vec{e}' = \vec{t}$ is the unit vector parallel to the tangent of the indicatrix. The equation $\vec{e}^*(s) = \vec{p}(s) \times \vec{e}(s)$ has infinity of solutions for the function $\vec{p}(s)$. If we take $\vec{p}_o(s)$ as a solution, the set of all solutions is given by $\vec{p}(s) = \vec{p}_o(s) + \lambda(s) \vec{e}(s)$, where $\lambda$ is a real scalar function of $s$. Therefore we have $\langle \vec{p}', \vec{e}' \rangle = \langle \vec{p}_o', \vec{e}' \rangle + \lambda$. By taking $\lambda = \lambda_o = -\langle \vec{p}_o', \vec{e}' \rangle$ we see that $\vec{p}_o(s) + \lambda_o(s) \vec{e}(s) = \vec{c}(s)$ is the unique solution for $\vec{p}(s)$ with $\langle \vec{c}', \vec{e}' \rangle = 0$. Then, the given dual curve $(\tilde{k})$ corresponding to the ruled surface

$$\varphi_e = \vec{c}(s) + v \vec{e}(s) \tag{12}$$

may be represented by

$$\tilde{e}(s) = \vec{e} + \varepsilon \vec{c} \times \vec{e} \tag{13}$$

where

$$\langle \vec{e}, \vec{e} \rangle = 1, \quad \langle \vec{e}', \vec{e}' \rangle = 1, \quad \langle \vec{c}', \vec{e}' \rangle = 0 . \tag{14}$$

Then we have

$$\|\tilde{e}'\| = \vec{t} + \varepsilon \det(\vec{c}', \vec{e}, \vec{t}) = 1 + \varepsilon \Delta \tag{15}$$

where $\Delta = \det(\vec{c}', \vec{e}, \vec{t})$. The dual arc-length $\bar{s}$ of the dual curve $(\tilde{k})$ is given by

$$\overline{s} = \int_0^s \|\tilde{e}'(u)\| du = \int_0^s (1+\varepsilon\Delta)du = s + \varepsilon\int_0^s \Delta du \qquad (16)$$

From (16) we have $\overline{s}' = 1 + \varepsilon\Delta$. Therefore, the dual unit tangent to the curve $\tilde{e}(s)$ is given by

$$\frac{d\tilde{e}}{d\overline{s}} = \frac{\tilde{e}'}{\overline{s}'} = \frac{\tilde{e}'}{1+\varepsilon\Delta} = \tilde{t} = \vec{t} + \varepsilon(\vec{c}\times\vec{t}) \qquad (17)$$

Introducing the dual unit vector $\tilde{g} = \tilde{e}\times\tilde{t} = \vec{g} + \varepsilon\vec{c}\times\vec{g}$ we have the dual frame $\{\tilde{e},\tilde{t},\tilde{g}\}$ which is known as dual geodesic trihedron or dual Darboux frame of $\varphi_e$ (or $(\tilde{e})$). Also, it is well known that the real orthonormal frame $\{\vec{e},\vec{t},\vec{g}\}$ is called the geodesic trihedron of the indicatrix $\vec{e}(s)$ with the derivations

$$\vec{e}' = \vec{t}, \quad \vec{t}' = \gamma\vec{g} - \vec{e}, \quad \vec{g}' = -\gamma\vec{t} \qquad (18)$$

where $\gamma$ is called the conical curvature[4]. Similar to (18), the derivatives of the vectors of the dual frame $\{\tilde{e},\tilde{t},\tilde{g}\}$ are given by

$$\frac{d\tilde{e}}{d\overline{s}} = \tilde{t}, \quad \frac{d\tilde{t}}{d\overline{s}} = \overline{\gamma}\tilde{g} - \tilde{e}, \quad \frac{d\tilde{g}}{d\overline{s}} = -\overline{\gamma}\tilde{t} \qquad (19)$$

where

$$\overline{\gamma} = \gamma + \varepsilon(\delta - \gamma\Delta), \quad \delta = \langle\vec{c}',\vec{e}\rangle \qquad (20)$$

and the dual darboux vector of the frame is $\tilde{d} = \overline{\gamma}\tilde{e} + \tilde{g}$. From the definition of $\Delta$ and (20) we also have

$$\vec{c}' = \delta\vec{e} + \Delta\vec{g} \qquad (21)$$

The dual curvature of dual curve(ruled surface) $\tilde{e}(s)$ is

$$\overline{R} = \frac{1}{\sqrt{(1+\overline{\gamma}^2)}} \qquad (22)$$

The unit vector $\tilde{d}_o$ with the same sense as the Darboux vector $\tilde{d} = \overline{\gamma}\tilde{e} + \tilde{g}$ is given by

$$\tilde{d}_o = \frac{\overline{\gamma}}{\sqrt{(1+\overline{\gamma}^2)}}\tilde{e} + \frac{1}{\sqrt{(1+\overline{\gamma}^2)}}\tilde{g} \qquad (23)$$

Then, the dual angle between $\tilde{d}_o$ and $\tilde{e}$ satisfies the followings

$$\cos\overline{\rho} = \frac{\overline{\gamma}}{\sqrt{(1+\overline{\gamma}^2)}}, \quad \sin\overline{\rho} = \frac{1}{\sqrt{(1+\overline{\gamma}^2)}} \qquad (24)$$

where $\overline{\rho}$ is the dual spherical radius of curvature. Hence $\overline{R} = \sin\overline{\rho}$, $\overline{\gamma} = \cot\overline{\rho}$ [9].

## 4. Darboux Approach to Bertrand Ruled Surfaces

Let $\varphi_e$ and $\varphi_{e_1}$ be two ruled surfaces generated by dual unit vectors $\tilde{e}$ and $\tilde{e}_1$ and let $\{\tilde{e}(\overline{s}), \tilde{t}(\overline{s}), \tilde{g}(\overline{s})\}$ and $\{\tilde{e}_1(\overline{s}_1), \tilde{t}_1(\overline{s}_1), \tilde{g}_1(\overline{s}_1)\}$ denote the dual Darboux frames of $\varphi_e$ and $\varphi_{e_1}$, respectively. Then $\varphi_e$ and $\varphi_{e_1}$ are called Bertrand surface offsets, if

$$\tilde{t}(\overline{s}) = \tilde{t}_1(\overline{s}_1) \qquad (25)$$

holds, where $\overline{s}$ and $\overline{s}_1$ are the dual arc-lengths of $\varphi_e$ and $\varphi_{e_1}$, respectively. By this definition, the relation between the trihedrons of the ruled surfaces $\varphi_e$ and $\varphi_{e_1}$ can be given as follows

$$\begin{pmatrix} \tilde{e}_1 \\ \tilde{t}_1 \\ \tilde{g}_1 \end{pmatrix} = \begin{pmatrix} \cos\bar{\theta} & 0 & -\sin\bar{\theta} \\ 0 & 1 & 0 \\ \sin\bar{\theta} & 0 & \cos\bar{\theta} \end{pmatrix} \begin{pmatrix} \tilde{e} \\ \tilde{t} \\ \tilde{g} \end{pmatrix} \quad (26)$$

where $\bar{\theta} = \theta + \varepsilon\theta^*$, $(0 \leq \theta \leq \pi,\ \theta^* \in \mathbb{R})$ is the dual angle between the generators $\tilde{e}$ and $\tilde{e}_1$ of Bertrand ruled surface $\varphi_e$ and $\varphi_{e_1}$. The angle $\theta$ is called the offset angle and $\theta^*$ is called the offset distance(Fig. 1).

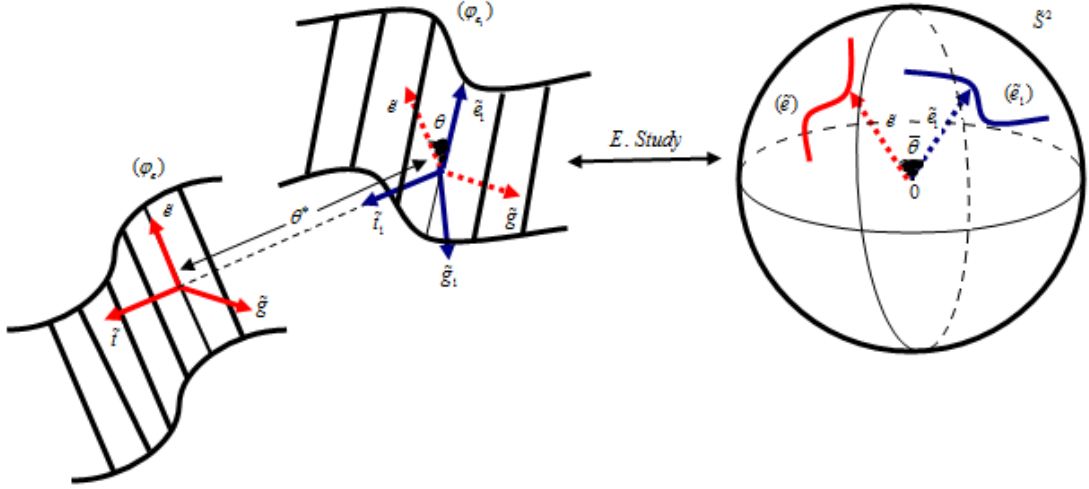

**Fig. 1.** E. Study Mapping of the Bertrand surface offsets.

Then, $\bar{\theta} = \theta + \varepsilon\theta^*$ is called dual offset angle of the Bertrand ruled surface $\varphi_e$ and $\varphi_{e_1}$. If $\theta = 0$ and $\theta = \pi/2$ then the Bertrand surface offsets are said to be oriented offsets and right offsets, respectively[5].

Now, we give some theorems and results characterizing Bertrand surface offsets. First we give the following theorem which is also given in [5] and proved by considering Blaschke frame. Then we give the new theorems and corollaries for Bertrand surface offsets.

**Theorem 4.1.** *Let $\varphi_e$ and $\varphi_{e_1}$ form a Bertrand surface offset. The offset angle $\theta$ and the offset distance $\theta^*$ are constants.*

**Proof.** Suppose that $\varphi_e$ and $\varphi_{e_1}$ form a Bertrand offset. Then from (26) we have

$$\tilde{e}_1 = \cos\bar{\theta}\,\tilde{e} - \sin\bar{\theta}\,\tilde{g}. \quad (27)$$

By differentiating (27) with respect to $\bar{s}$ we get

$$\frac{d\tilde{e}_1}{d\bar{s}} = (\cos\bar{\theta} + \bar{\gamma}\sin\bar{\theta})\tilde{t}_1 - \frac{d\bar{\theta}}{d\bar{s}}\tilde{g}_1. \quad (28)$$

Since $\dfrac{d\tilde{e}_1}{d\bar{s}}$ and $\tilde{g}_1$ are orthogonal, from (28) we have $\dfrac{d\bar{\theta}}{d\bar{s}} = 0$, i.e., $\bar{\theta} = \theta + \varepsilon\theta^*$ is constant. Then, we obtain that the offset angle $\theta$ and the offset distance $\theta^*$ are constants.

**Theorem 4.2.** *If $\varphi_e$ and $\varphi_{e_1}$ form a Bertrand surface offset, then there is the following differential relationship between the dual arc-length parameters of $\varphi_e$ and $\varphi_{e_1}$*

$$\frac{d\bar{s}_1}{d\bar{s}} = \cos\bar{\theta} + \bar{\gamma}\sin\bar{\theta}. \tag{29}$$

**Proof.** Suppose that $\varphi_e$ and $\varphi_{e_1}$ form a Bertrand offset. Then, from (26) and Theorem 4.1 we have

$$\frac{d\tilde{e}_1}{d\bar{s}_1} = (\cos\bar{\theta} + \bar{\gamma}\sin\bar{\theta})\frac{d\bar{s}}{d\bar{s}_1}\tilde{t}_1. \tag{30}$$

Since $\frac{d\tilde{e}_1}{d\bar{s}_1} = \tilde{t}_1$, from (30) we obtain Eq. (29).

Real part of (29) is obtain as follows

$$\frac{ds_1}{ds} = \cos\theta + \gamma\sin\theta, \tag{31}$$

and it gives the relation between the arc-length of the real indicatrixes $\vec{e}(s)$ and $\vec{e}_1(s_1)$.

**Theorem 4.3.** *Let $\varphi_e$ and $\varphi_{e_1}$ form a Bertrand surface offset. Then there are the following relationships between the arc-length parameters of $\varphi_e$ and $\varphi_{e_1}$*

$$\begin{cases} s_1 = \cos\theta\, s + \left(\int_0^{\bar{s}}\gamma ds\right)\sin\theta \\ \int_0^{s_1}\Delta_1 du_1 = \left(\int_0^{\bar{s}}(\gamma\Delta+\gamma^*)ds\right)\sin\theta + \left(\int_0^s\Delta du\right)\cos\theta + \theta^*\left[\left(\int_0^{\bar{s}}\gamma ds\right)\cos\theta - \sin\theta\, s\right] \end{cases} \tag{32}$$

**Proof.** If $\varphi_e$ and $\varphi_{e_1}$ form a Bertrand surface offset, then Theorem 4.2 holds and from (29) we have

$$\bar{s}_1 = \int_0^{\bar{s}}\cos\bar{\theta}\, d\bar{s} + \int_0^{\bar{s}}\sin\bar{\theta}\,\bar{\gamma}d\bar{s}. \tag{33}$$

Since $\bar{\theta}$ is constant, equality (33) may be written as

$$\bar{s}_1 = \bar{s}\cos\bar{\theta} + \sin\bar{\theta}\int_0^{\bar{s}}\bar{\gamma}d\bar{s}. \tag{34}$$

By separating the right side of (34) into real and dual parts and using the fact that $s^* = \int_0^s \Delta du$ we obtain

$$\bar{s}_1 = s\cos\theta + \left(\int_0^{\bar{s}}\gamma ds\right)\sin\theta + \varepsilon\left\{\left(\int_0^{\bar{s}}(\gamma\Delta+\gamma^*)ds\right)\sin\theta + \left(\int_0^s\Delta du\right)\cos\theta + \theta^*\left[\left(\int_0^{\bar{s}}\gamma ds\right)\cos\theta - s\sin\theta\right]\right\}. \tag{35}$$

From (16) we have $\bar{s}_1 = s_1 + \varepsilon\int_0^{s_1}\Delta_1 du_1$. Then from (35) we find the equalities (32).

**Corollary 4.1.** *Let $\varphi_e$ and $\varphi_{e_1}$ form a Bertrand surface offset and let $\varphi_e$ be a developable ruled surface generated by dual unit vector $\tilde{e}$. Then, $\varphi_{e_1}$ is developable if and only if*

$$\left(\int_0^{\bar{s}}\gamma^* ds\right)\sin\theta + \theta^*\left[\left(\int_0^{\bar{s}}\gamma ds\right)\cos\theta - \sin\theta\, s\right] = 0 \tag{36}$$

*holds.*

**Proof.** Assume that $\varphi_e$ and $\varphi_{e_1}$ form a Bertrand surface offset and let $\varphi_e$ be developable. Then we have $\Delta = 0$. From Theorem 4.3, the condition to satisfy $\Delta_1 = 0$ is that

$$\left(\int_0^{\bar{s}} \gamma^* ds\right) \sin\theta + \theta^* \left[\left(\int_0^{\bar{s}} \gamma ds\right) \cos\theta - \sin\theta s\right] = 0.$$

Conversely if $\varphi_e$ is developable and (36) holds, then from (32) we have that $\Delta_1 = 0$, i.e, $\varphi_{e_1}$ is developable.

From (32) we can give the following result

**Result 4.1. i)** *If the ruled surfaces $\varphi_e$ and $\varphi_{e_1}$ form an oriented offset then the real indicatrixes $\vec{e}(s)$ and $\vec{e}_1(s_1)$ have the same arc length parameter, i.e, $s = s_1$ and*

$$\int_0^{s_1} \Delta_1 du_1 = \int_0^s \Delta du + \theta^* \int_0^{\bar{s}} \gamma ds.$$

**ii)** *If the ruled surfaces $\varphi_e$ and $\varphi_{e_1}$ form a right offset then the following equalities holds*

$$s_1 = \int_0^{\bar{s}} \gamma ds, \quad \int_0^{s_1} \Delta_1 du_1 = \int_0^{\bar{s}} (\gamma\Delta + \gamma^*) ds - \theta^* s.$$

**Theorem 4.4.** *Let $\varphi_e$ and $\varphi_{e_1}$ form a Bertrand surface offset. The offset distance $\theta^*$ is found as follows*

$$\theta^* = \frac{\delta \cos\theta - \Delta \sin\theta}{\cos\theta + \gamma \sin\theta} - \delta_1. \tag{37}$$

**Proof.** Let $\varphi_e$ and $\varphi_{e_1}$ form a Bertrand surface offset and the striction lines of $\varphi_e$ and $\varphi_{e_1}$ be $c(s)$ and $c_1(s_1)$, respectively. Then by the definition of Bertrand offset we can write

$$\vec{c}_1 = \vec{c} + \theta^* \vec{t}. \tag{38}$$

Differentiating (38) with respect to $s_1$ we have

$$\frac{d\vec{c}_1}{ds_1} = \left(\frac{d\vec{c}}{ds} + \theta^* (\gamma\vec{g} - \vec{e})\right) \frac{ds}{ds_1}. \tag{39}$$

From (20) we have $\delta_1 = \langle d\vec{c}_1/ds_1, e_1 \rangle$. Then from the real part of (26) and (39) we obtain

$$\delta_1 = \left(\cos\theta \langle d\vec{c}/ds, \vec{e} \rangle - \theta^* \cos\theta - \sin\theta \langle d\vec{c}/ds, \vec{g} \rangle - \theta^* \gamma \sin\theta\right) \frac{ds}{ds_1}. \tag{40}$$

Using the equalities

$$\delta = \langle d\vec{c}/ds, e \rangle, \quad \langle d\vec{c}/ds, \vec{g} \rangle = \langle d\vec{c}/ds, \vec{e} \times \vec{t} \rangle = \det(d\vec{c}/ds, \vec{e}, \vec{t}) = \Delta, \tag{41}$$

from (40) we have

$$\theta^* = \frac{\delta \cos\theta - \Delta \sin\theta}{\cos\theta + \gamma \sin\theta} - \delta_1.$$

**Result 4.2. i)** *If the ruled surfaces $\varphi_e$ and $\varphi_{e_1}$ form an oriented offset then $\theta^* = \delta - \delta_1$ holds.*

**ii)** *If the ruled surfaces $\varphi_e$ and $\varphi_{e_1}$ form a right offset then $\theta^* = \frac{\Delta}{\gamma} - \delta_1$ holds.*

**Theorem 4.5.** *Let $\varphi_e$ and $\varphi_{e_1}$ form a Bertrand surface offset. The relationship between $\Delta$ and $\Delta_1$ is given by*

$$\Delta_1 = \frac{\Delta\cos\theta + \delta\sin\theta + \theta^*(\gamma\cos\theta - \sin\theta)}{\cos\theta + \gamma\sin\theta}. \qquad (42)$$

**Proof.** We know that $\Delta_1 = \det(d\vec{c}_1/ds_1, \vec{e}_1, \vec{t}_1)$. Then from (26) and (39) we obtain

$$\begin{aligned}\Delta_1 &= \frac{ds}{ds_1}\det\left(\left(\frac{d\vec{c}}{ds} + \theta^*(\gamma\vec{g} - \vec{e})\right), \cos\theta\vec{e} - \sin\theta\vec{g}, \vec{t}\right) \\ &= \frac{ds}{ds_1}\left[\det\left(\frac{d\vec{c}}{ds}, \vec{e}, \vec{t}\right)\cos\theta - \det\left(\frac{d\vec{c}}{ds}, \vec{g}, \vec{t}\right)\sin\theta + \theta^*\left(\det(\vec{g}, \vec{e}, \vec{t})\gamma\cos\theta + \det(\vec{e}, \vec{g}, \vec{t})\sin\theta\right)\right].\end{aligned} \qquad (43)$$

Writing $\Delta = \det(d\vec{c}/ds, \vec{e}, \vec{t})$ and

$$\det(\vec{e}, \vec{g}, \vec{t}) = -\det(\vec{g}, \vec{e}, \vec{t}) = -\langle \vec{g}, \vec{e}\times\vec{t}\rangle = -1,$$

$$\det\left(\frac{d\vec{c}}{ds}, \vec{g}, \vec{t}\right) = \left\langle\frac{d\vec{c}}{ds}, \vec{g}\times\vec{t}\right\rangle = -\left\langle\frac{d\vec{c}}{ds}, \vec{e}\right\rangle = -\delta.$$

in (43), we have equality (42).

**Corollary 4.2.** *Let $\varphi_e$ and $\varphi_{e_1}$ form a Bertrand surface offset and let $\varphi_e$ be developable. Then, the Bertrand offset $\varphi_{e_1}$ is developable if and only if*

$$\theta^* = \frac{\delta\sin\theta}{\sin\theta - \gamma\cos\theta}, \qquad (44)$$

*holds.*

**Proof.** Let $\varphi_e$ be a developable ruled surface, i.e, $\Delta = 0$. Then from Theorem 4.5, the condition to satisfy $\Delta_1 = 0$ is that

$$\theta^* = \frac{\delta\sin\theta}{\sin\theta - \gamma\cos\theta}.$$

Conversely if $\varphi_e$ is developable and (44) holds, then from (42) we have that $\Delta_1 = 0$, i.e, $\varphi_{e_1}$ is developable.

From Theorem 4.5 we can give the following result.

**Result 4.3. i)** *If the ruled surfaces $\varphi_e$ and $\varphi_{e_1}$ form an oriented offset, then $\Delta_1 = \Delta + \theta^*\gamma$ holds.*

**ii)** *If the ruled surfaces $\varphi_e$ and $\varphi_{e_1}$ form a right offset, then $\Delta_1 = \dfrac{\delta - \theta^*}{\gamma}$ holds.*

**Theorem 4.6.** *Let $\varphi_e$ and $\varphi_{e_1}$ form a Bertrand surface offset. There is the following relationship between the conical curvatures $\gamma$ and $\gamma_1$,*

$$\gamma_1 = \frac{\gamma\cos\theta - \sin\theta}{\cos\theta + \gamma\sin\theta}. \qquad (45)$$

**Proof.** From (18) and (26) we have

$$\begin{aligned}
\gamma_1 &= -\langle \vec{g}_1', \vec{t}_1 \rangle \\
&= -\left\langle \frac{d}{ds_1}(\sin\theta \vec{e} + \cos\theta \vec{g}), \vec{t} \right\rangle \\
&= (\gamma\cos\theta - \sin\theta)\frac{ds}{ds_1}.
\end{aligned}$$

From (31) we have (45).

**Corollary 4.3.** *Let $\varphi_e$ and $\varphi_{e_1}$ form a Bertrand surface offset. Then the offset angle is given by*

$$\theta = \arctan\left(\frac{\gamma - \gamma_1}{1 + \gamma\gamma_1}\right). \tag{46}$$

**Proof.** From (44) we have (46) immediately.

**Result 4.4. i)** *If the ruled surfaces $\varphi_e$ and $\varphi_{e_1}$ form an oriented offset, then $\gamma_1 = \gamma$.*

**ii)** *If the ruled surfaces $\varphi_e$ and $\varphi_{e_1}$ form a right offset, then $\gamma\gamma_1 = -1$.*

**Theorem 4.7.** *If $\varphi_e$ and $\varphi_{e_1}$ form a Bertrand surface offset, then there exists the following relationship between $\overline{\gamma}$ and $\overline{\gamma}_1$,*

$$\begin{aligned}
\overline{\gamma}_1 = \frac{ds}{ds_1}(\cos\theta\gamma - \sin\theta) + \frac{ds}{ds_1}\varepsilon\bigg[&\left(\frac{ds}{ds_1}\cos\theta(\sin\theta - \gamma\cos\theta) - \sin\theta\right)\Delta \\
&+\left(\cos\theta - \frac{ds}{ds_1}\sin\theta(\gamma\cos\theta - \sin\theta)\right)\delta \\
&-\left(\frac{ds}{ds_1}(\gamma\cos\theta - \sin\theta)^2 + \cos\theta + \gamma\sin\theta\right)\theta^*\bigg].
\end{aligned} \tag{47}$$

**Proof.** From (20), (37), (42) and (45), by direct calculation we have (47)

**Result 4.5.i)** *If the ruled surfaces $\varphi_e$ and $\varphi_{e_1}$ form an oriented offset then $\overline{\gamma}_1 = \overline{\gamma} + \varepsilon\left(-\gamma\Delta + \delta - (1+\gamma^2)\theta^*\right)$ holds.*

**ii)** *If the ruled surfaces $\varphi_e$ and $\varphi_{e_1}$ form a right offset then*

$$\overline{\gamma}_1 = \frac{1}{\gamma}\left[-1 + \varepsilon\left(-\Delta + \frac{1}{\gamma}\delta - \frac{\gamma^2+1}{\gamma}\theta^*\right)\right], \tag{48}$$

*holds.*

**Theorem 4.8.** *Let $\varphi_e$ and $\varphi_{e_1}$ form a Bertrand surface offset. Then, the dual curvature of $\varphi_{e_1}$ is given by*

$$\overline{R}_1 = \left[1 + A^2\right]^{-1/2} - B\varepsilon\left[1 + A^2\right]^{-3/2}, \tag{49}$$

*where*

$$\begin{cases} A = \dfrac{ds}{ds_1}(\gamma\cos\theta - \sin\theta), \\ \\ B = A\dfrac{ds}{ds_1}\big[-(\sin\theta + A\cos\theta)\Delta + (\cos\theta - A\sin\theta)\delta \\ \qquad - (A(\gamma\cos\theta + \sin\theta) + \cos\theta + \gamma\sin\theta)\theta^*\big]. \end{cases} \quad (50)$$

**Proof.** From (22) and (48), (49) and (50) are obtained immediately.

**Theorem 4.9.** *Let $\varphi_e$ and $\varphi_{e_1}$ form a Bertrand surface offset. Then, the dual spherical radius of curvature of $\varphi_{e_1}$ is given by*

$$\cos\bar{\rho}_1 = A(1+A^2)^{-1/2} + \varepsilon BA^{-1}(1+A^2)^{-3/2}, \quad (51)$$

*where $A$ and $B$ are same as in (50).*

**Proof.** From (24) we have $\sin\bar{\rho} = \bar{R}$. Then by (4) and (47) we get

$$\sin\rho_1 = \left[1+A^2\right]^{-1/2}, \quad \rho_1^*\cos\rho_1 = -B\left[1+A^2\right]^{-3/2}. \quad (52)$$

By using the trigonometric relation $\sin^2\rho_1 + \cos^2\rho_1 = 1$ and considering (52), we have

$$\cos\rho_1 = A\left[1+A^2\right]^{-1/2}. \quad (53)$$

By writing (53) in the second equality of (52) we obtain

$$\rho_1^* = \dfrac{-B}{A(1+A^2)}. \quad (54)$$

Then by using the extension $\cos\bar{\rho}_1 = \cos\rho_1 - \varepsilon\rho_1^*\sin\rho_1$, from (52)-(54) we have (51).

**Corollary 4.4.** *The relationship between the offset angle $\theta$ and the real spherical radius of curvature $\rho_1$ is given by*

$$\cot\rho_1 = \dfrac{\gamma\cos\theta - \sin\theta}{\cos\theta + \gamma\sin\theta}. \quad (55)$$

**Proof.** The proof is clear from (31), (52) and (53).

**Result 4.6. i)** *If the ruled surfaces $\varphi_e$ and $\varphi_{e_1}$ form an oriented offset then the relationship between the real spherical radius of curvature $\rho_1$ and conical curvature $\gamma$ is given by*

$$\cot\rho_1 = \gamma.$$

**ii)** *If the ruled surfaces $\varphi_e$ and $\varphi_{e_1}$ form a right offset then the relationship between the real spherical radius of curvature $\rho_1$ and conical curvature $\gamma$ is given by*

$$\cot\rho_1 = \dfrac{1}{\gamma}.$$

## 5. Conclusions

In this paper, we give the characterizations of Bertrand offsets of ruled surfaces in view of dual geodesic trihedron(dual Darboux frame). We find new relations between the invariants of Bertrand surface offsets. Furthermore, we give the relationships for Bertrand surface offsets to be developable.